\documentclass[12pt]{article}
\usepackage{amsmath,amsthm,amssymb}

\numberwithin{equation}{section} \theoremstyle{remark}

\numberwithin{equation}{section}

\def\1{\bf 1}

\usepackage[dvips]{graphicx}
\usepackage{amssymb}
\usepackage{amsmath}
\usepackage{latexsym}
\usepackage[cp1251]{inputenc}
\usepackage[english]{babel}
\begin{document}

{\bf S. Albeverio$^{1,2,3,4}$, V. Koshmanenko$^5$, and I.
Samoilenko$^6$}

\begin{center}
\large \textbf{The conflict interaction between two complex
systems. Cyclic migration}\\

\ \ \\ \ \ \ \

\end{center}

\begin{abstract}
We construct and study a discrete time model describing the
conflict interaction between two complex systems with non-trivial
internal structures.  The external conflict interaction is based
on the model of alternative interaction between a pair of
non-annihilating opponents. The internal conflict dynamics is
similar to the one of a predator-prey model. We show that the
typical trajectory of the complex system converges to an
asymptotic attractive cycle. We propose an interpretation of our
model in terms of migration processes.
\end{abstract}

$^1$ {Institut f\"{u}r Angewandte Mathematik, Universit\"{a}t
Bonn, Wegelerstr. 6, D-53115 Bonn\ (Germany); }$^2${SFB 611, \
Bonn, \ BiBoS, (Bielefeld - Bonn); }$^3${IZKS Bonn}; $^4${CERFIM,
Locarno and Acc. Arch. (USI) (Switzerland)} {e-mail:
albeverio@uni-bonn.de}

$^5${ Institute of Mathematics, Tereshchenkivs'ka str. 3, Kyiv
01601 Ukraine} \ {e-mail: kosh@imath.kiev.ua}

$^6${Institute of Mathematics, Tereshchenkivs'ka str. 3, Kyiv
01601 Ukraine} \ {e-mail: isamoil@imath.kiev.ua}

\textbf{2000 Mathematics Subject Classification: 91A05,  91A10,
90A15,  90D05,  37L30,  28A80}

\textbf{Key words: } Lotka-Volterra equations, predator-prey
model, conflict interaction, dynamical system, cyclic attractor,
 limiting distributions, migration

\section{Introduction}
Since the beginning of 20-th century the Lotka-Volterra model of
prey-predator interaction is one of the main models for simulation
of many processes in population theory and economics. As a rule,
continuous models where Lotka-Volterra equations have
ratio-depended parameters are studied (see, for example \cite{BC,
CM, CKG, Ku, KuBe, LSGM, Mur, StOl, Tu}). Logistical and Ricker's
models are also studied in some works, for example \cite{CM}. In
the majority of works the prey-predator interaction is treated
only inside of single region, and no migration from one region to
another is considered.

In some works \cite{CKG, CM} models with migration are studied,
with a process of migration.

There are also only few works (see \cite{CM} and references
wherein), in which discrete models are considered, though in
reality such processes are more natural, since they take better
into account seasonal phenomena (reproduction, migration, etc.).

The main aim of the majority of works is determination of stable
points, bifurcation points, asymptotic behavior, and analysis of
model's depending on the coefficients of the equations.

In \cite{BC} synchronization of population dynamics with natural
phenomena (like change of seasons and floods) is studied. In the
work \cite{Tu} dependence of population dynamics on population
density and migration is studied. In the work \cite{CKG} migration
is not assumed to be random, but aims at maximization of some
function of the population. At last, in \cite{BC} the influence of
stochastic terms in a Lotka-Volterra model is described, and
interesting figures are presented.

In recent works \cite{SaTa1, SaTa} Salam and Takahashi study
conflict models, similar to ours. They introduce important and more
complex multi-opponent systems. In \cite{SaTa1} not only conflict,
but cooperation between opponents is studied. The figures, obtained
by them, are very similar to Figure 3 of the present work.

In this work we construct a model that joins two most rarely
studied variants of Lotka-Volterra model, i.e., a discrete model
with migration. Here individuals migrate not randomly, but
according to strategies, discussed in section 5.

We construct the model of the conflict interaction between a pair
of complex systems A and B. The system is a finite set of positive
numbers: $\mathbf{P}=(P_1,\ldots, P_N)$ for A and
$\mathbf{R}=(R_1,\ldots,R_N)$ for B, where $N$ means the quantity
of parameters that characterize the system. We study dynamics in
the discrete time. So, the evolution of every system is described
by the sequence of vectors with non-negative coordinates
$\mathbf{P}^n=(P^{(n)}_1,\ldots, P^{(n)}_N)$ for A, and
$\mathbf{R}^n=(R^{(n)}_1,\ldots, R^{(n)}_N)$ for B,
$n=1,2,\ldots$. The vectors $\mathbf{P}$ and $\mathbf{R}$
correspond to the moment $n=0$. Naturally, each system tries to
reach the optimal values of its coordinates. In reality, due to
the conflict interaction, every coordinate changes in a
complicated way. The evolution of all changes is determined by
double dependence: by the conflict interaction between systems
(which we shortly describe in section 3), and by the mutual
"fight" of coordinates (of the prey-predator type interaction)
inside every system. We suppose that every system is complex in
the sense that its elements may be treated as one of the types:
dominant (predators, employers) or dependent (preys, workers). So,
every coordinate $P^{(n)}_i,R^{(n)}_i$ may be regarded as the
quantity (population) of dominant, respectively dependent species
at the position $i$ at time $n$.

The law of evolution inside of each (independent) system is
described in section 2. We suppose this law is identical in every
system and is based on the well-known discrete Lotka-Volterra
equation.

In section 3 we shortly call the main results on conflict
interaction between non-annihilating opponents.

In section 4, that includes the main results of the work, we
construct a dynamical system describing simultaneous conflict
interaction both inside every system and between the systems. The
outer interaction is an alternative conflict between
non-annihilating systems, whereas the inner one is a prey-predator
model of Lotka-Volterra type.

We may join these two types of interactions in a discrete time.
Thus, our dynamical system consists of a discrete sequence of
states. Two operations happen at any fixed moment of time:
redistribution of probabilities to occupy some controversial
positions by opponent systems, and quantitative changes (namely
population) of all species inside both systems.

The computer modelling of such a complex interaction shows some
very interesting phenomena. In this work we limit however
ourselves to present only one observation. Namely, under an
appropriate choice of parameters and initial data the complex
system oscillates. We find a rather wide range for initial data
for which the population trajectory in phase-space becomes cyclic.
Moreover, we observe the stability of the limit cycle, so it is an
attractor.

\section{Traditional models of population dynamics} Malthus proposed in 1798 the population equation
$$\frac{dP}{dt}=(b-d)P, \eqno(2.1)$$ where $P$
is the cumulative number of individuals (species), and $b,d$ are
the natural birth and death rates. In reality, one expects
exponential solution
$$P(t)=P(0)e^{(b-d)t}$$ describes the ideal population of
biological species. The exponential rise, if $b-d>0$, or decrease,
if $b-d<0$ at most in a local period of time.

Verhulst introduced in 1838 a more realistic equation with
saturation terms: $$\frac{dP}{dt}=(b-d)P-cP^2, \eqno(2.2)$$ where
the coefficient $c>0$ represents the competition activity of
individuals for living resources. The square power corresponds a
conception of an alternative law of access to the living resource.

The solution of (2.2) describes the $S$-shaped logistical curve
(see Figure  1) and corresponds better to the actual behavior of
many population processes.

\begin{figure}
\centering
\includegraphics{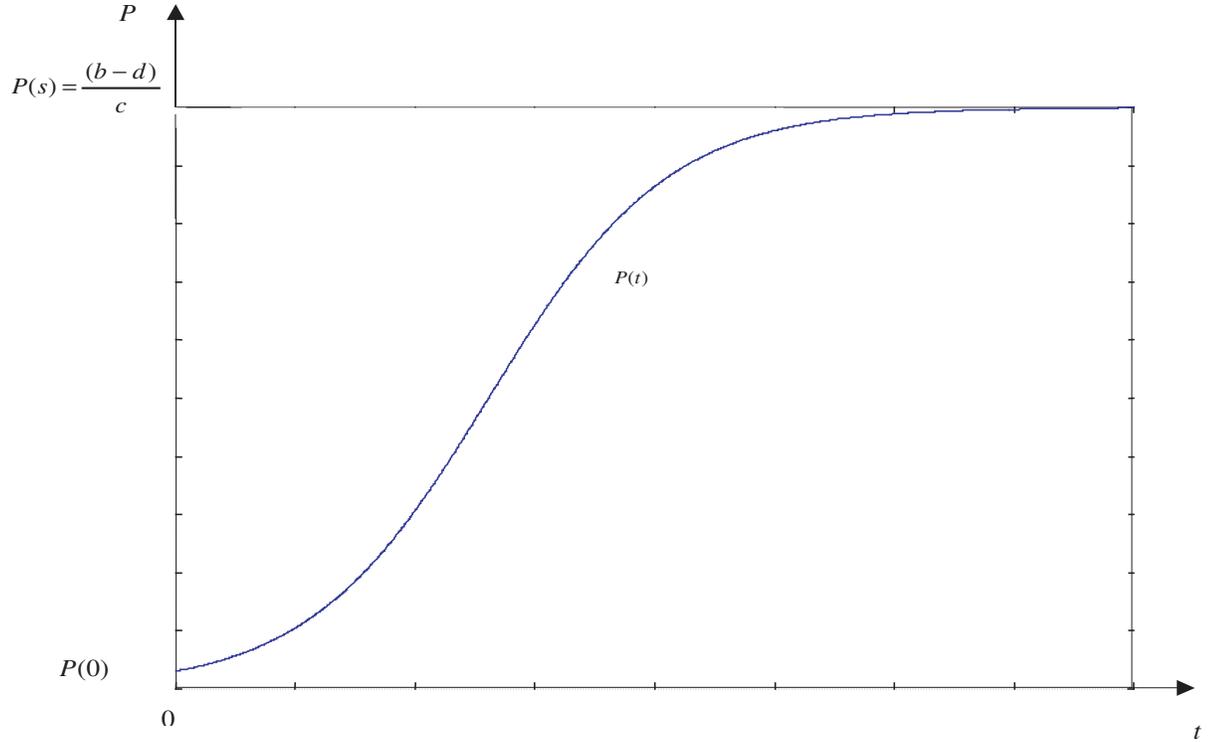}
\caption{{\it A typical shape of the logistical curve
$\frac{dP}{dt}=(b-d)P-cP^2.$ }}
\end{figure}

The curve starts with a small value $P(0)$, exponentially
increases, and then saturate at the capacity $P(s)=\frac{b-d}{c},
b-d>0.$

In the economic context, equation (2.2) can be written as follows
$$\frac{dM}{dt}=(g-l)M-fM^2, \eqno(2.3)$$ where $M$ is the capital
(money), $g$ and $l$ are the average gain and loss percentages on
the capital, and $f$ stands for the coefficient of confrontation
between individuals. If $g-l<0$, the capital decays to 0
exponentially; if $g-l>0,$ at the beginning the capital increases
exponentially quick, but then the growth slows down, so that it
never reaches asymptotic value of saturation $M(s)=\frac{g-l}{f}$.

Lotka (1907) and Volterra (1901) extended the Verhulst logistical
equation to the Lotka-Volterra equations intended for the
description of amount changes in populations of two species in
interaction. These equations are also known under the name of
predator-prey model. We will refer on Lotka-Volterra equations in
the following form:
$$\begin{array}{c}
  \frac{dP}{dt}=aP-bPR-cP^2 \\
  \frac{dR}{dt}=-dR+ePR-fR^2,
\end{array} \eqno(2.4)$$ where all coefficients are nonnegative.

The population of prey is described by the first equation. Without
presence of any predators it grows exponentially at the beginning
and then comes to the fixed capacity $P(s)=a/c.$ The predators,
without any prey to feed on, die out. When both species are
present, the growth of the prey is limited by the predators, due
to the term $-bPR$, and the predators grow if the amount of prey
available, i.e. if $ePR$ is large enough.

There are many publications devoted to the analysis of
Lotka-Volterra equations (2.4) (see for example \cite{Mur} and
references wherein).

The models with discrete time are also studied. In this case,
equations (2.4) have the following view: $$\begin{array}{c}
  P_1^{(n)}=P_1^{(n-1)}+P_1^{(n-1)}(a-bP_2^{(n-1)}-cP_1^{(n-1)}), \\
  P_2^{(n)}=P_2^{(n-1)}+P_2^{(n-1)}(-d+eP_1^{(n-1)}-fP_2^{(n-1)}).
\end{array} \eqno(2.5)$$

Typical behaviour of discrete Lotka-Volterra model is shown in
Figure 2.

\begin{figure}
\centering
\includegraphics{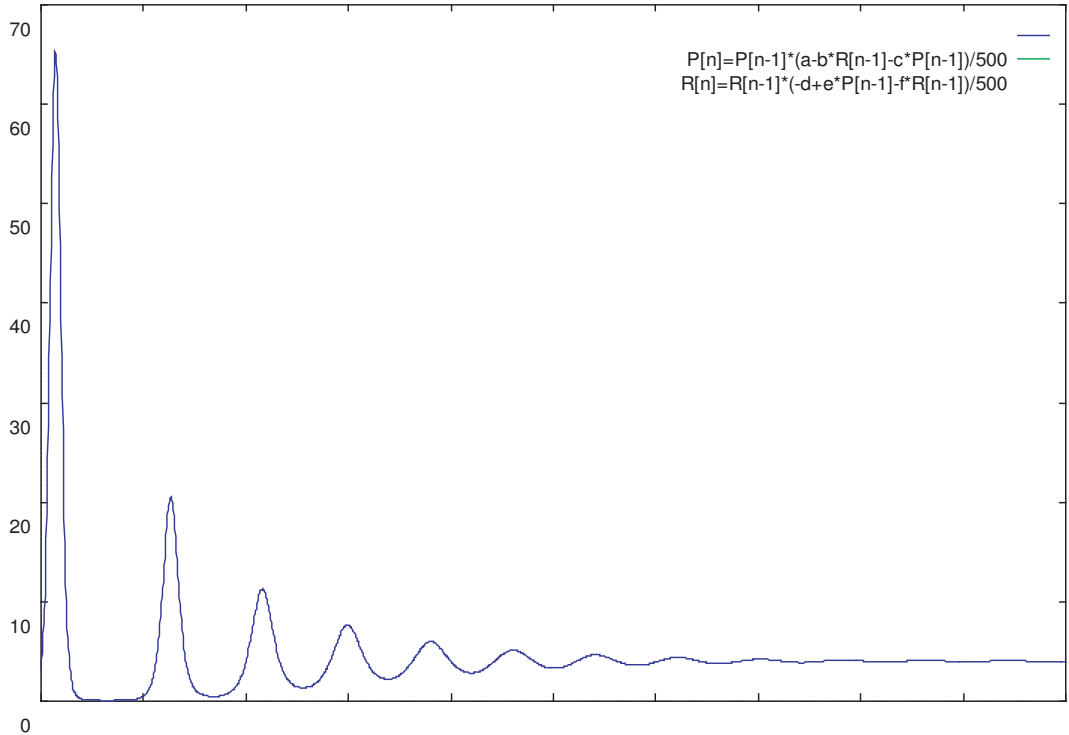}
\caption{{\it Lotka-Volterra model with discrete time}}
$P_1^{(n)}=P_1^{(n-1)}(a-bP_1^{(n-1)}-cP_1^{(n-1)})$
$$P_2^{(n)}=P_2^{(n-1)}(-d+eP_2^{(n-1)}-fP_2^{(n-1)})$$
$$a=0.2,b=0.006,c=0.002,d=0.008,e=0.002,f=0,$$
$$P_1^{(0)}=3,P_2^{(0)}=5.$$
\end{figure}

\section{Conflict interaction between non-annihilating opponents}
In this section we shortly remind an alternative approach to
describe the redistribution of conflicting positions between two
opponents, say A and B, concerning an area of common interests.

We consider the simplest case where the existence space of common
interests is a finite set of positions
$\Omega=\{\omega_1,\ldots,\omega_N\}, N\geq 2.$ Each of the
opponents A and B tries to occupy a position $\omega_i,
i=1,\ldots,N$ with a probability $P_A(\omega_i)=p_i\geq 0$ or
$P_B(\omega_i)=r_i\geq 0.$ The starting distributions of A and B
along $\Omega$ are arbitrary and normed:
$\sum_{i=1}^Np_i=1=\sum_{i=1}^Nr_i.$ A and B can not be present
simultaneously in a same position $\omega_i$. The interaction
between A and B is considered in discrete time $t\in
\mathbb{N}_0.$ We introduce the noncommutative conflict
composition between stochastic vectors
$\mathbf{p}^0=(p_1,\ldots,p_N), \mathbf{r}^0=(r_1,\ldots,r_N) \in
\mathbb{R}_+^N$:
$$\mathbf{p}^1:=\mathbf{p}^0*\mathbf{r}^0,
\mathbf{r}^1=\mathbf{r}^0*\mathbf{p}^0,
\mathbf{p}^0\equiv\mathbf{p}, \mathbf{r}^0\equiv\mathbf{r},$$
where the coordinates of $\mathbf{p}^1, \mathbf{r}^1$ are defined
as follows
\begin{equation}\label{result_Aa}
p_i^{(1)}=\frac{p_i^{(0)}(1-\alpha
r_i^{(0)})}{1-\alpha\sum_{i=1}^Np_i^{(0)}r_i^{(0)}},
r_i^{(1)}=\frac{r_i^{(0)}(1-\alpha
p_i^{(0)})}{1-\alpha\sum_{i=1}^Np_i^{(0)}r_i^{(0)}},
\end{equation} where the coefficient $-1\leq\alpha\leq 1$,
$\alpha\neq0$ stands for the activity interaction. At the $n$th
step of the conflict dynamics we get two vectors
$$\mathbf{p}^n=\mathbf{p}^{n-1}*\mathbf{r}^{n-1}\equiv
\mathbf{p}^0*^n\mathbf{r}^0,
\mathbf{r}^n=\mathbf{r}^{n-1}*\mathbf{p}^{n-1}\equiv
\mathbf{r}^0*^n\mathbf{p}^0$$ with coordinates
$$p_i^{(n)}=\frac{p_i^{(n-1)}(1-\alpha r_i^{(n-1)})}{z_n},
r_i^{(n)}=\frac{r_i^{(n-1)}(1-\alpha p_i^{(n-1)})}{z_n},$$ with
$z_n$ a normalization coefficient given by
$$z_n=1-\alpha(\mathbf{p}^{n-1},\mathbf{r}^{n-1}),$$ with
$(\cdot,\cdot)$ the inner product in $\mathbb{R}^N.$

The behavior of the state $\{\mathbf{p}^n,\mathbf{r}^n\}$ at time
$t=n$ for $n\to\infty$ has been investigated in \cite{ABodK,
BKoDo, KoTC, KoTC1, KoDo, KoKh}. We shortly describe the results.

{\bf Theorem.1.} {\it For any pair of non-orthogonal stochastic
vectors ${\bf p},{\bf r} \in \mathbb{R}_+^N,$ $({\bf p},{\bf r})>
0$, and fixed interaction intensivity parameter $\alpha\neq0$,
$-1\leq\alpha\leq 1,$ with condition $\alpha\neq \frac1{({\bf
p},{\bf r})}$, the sequence of states
$\{\mathbf{p}^n,\mathbf{r}^n\}$ tends to the limit state $\{{\bf
p}^{\infty},{\bf r}^{\infty}\}$ $${\bf p}^{\infty}= \lim_{n \to
\infty}{\bf p}^{n}, \ \ {\bf r}^{\infty}= \lim_{n \to \infty}{\bf
r}^{n}.$$ This limit state is invariant with respect to the
conflict interaction: $$ {\bf p}^{\infty}= {\bf p}^{\infty} {
\divideontimes} {\bf r}^{\infty}, \ \ {\bf r}^{\infty}= {\bf
r}^{\infty} { \divideontimes} {\bf p}^{\infty}. $$ Moreover, $$
\left\{
\begin{array}{cccc}{\bf p}^{\infty} \perp
{\bf r}^{\infty}, & \textrm{if \  \ ${\bf p} \neq {\bf r}$} \ \
\textrm{ and \ \  $0<\alpha\leq 1$} \\
 {\bf p}^{\infty}=
{\bf r}^{\infty}, & \textrm{ in all other cases.}
\end{array}\\
  \right.
$$ }

We emphasize that in the case of a purely repulsive interaction,
$0<\alpha\leq 1,$ if the starting distributions are different, then
the limiting vectors are orthogonal. Therefore each of the vectors
$\mathbf{p}^{\infty}, \mathbf{r}^{\infty}$ contains by necessity
some amount of zero coordinates on different positions $\omega_i$.
For example the typical limiting picture for $\mathbf{p}^n,
\mathbf{r}^n\in\mathbb{R}^3_+$ is presented in Figure 3 (comp. with
\cite{SaTa1, SaTa}).

\begin{figure}
\centering
\includegraphics{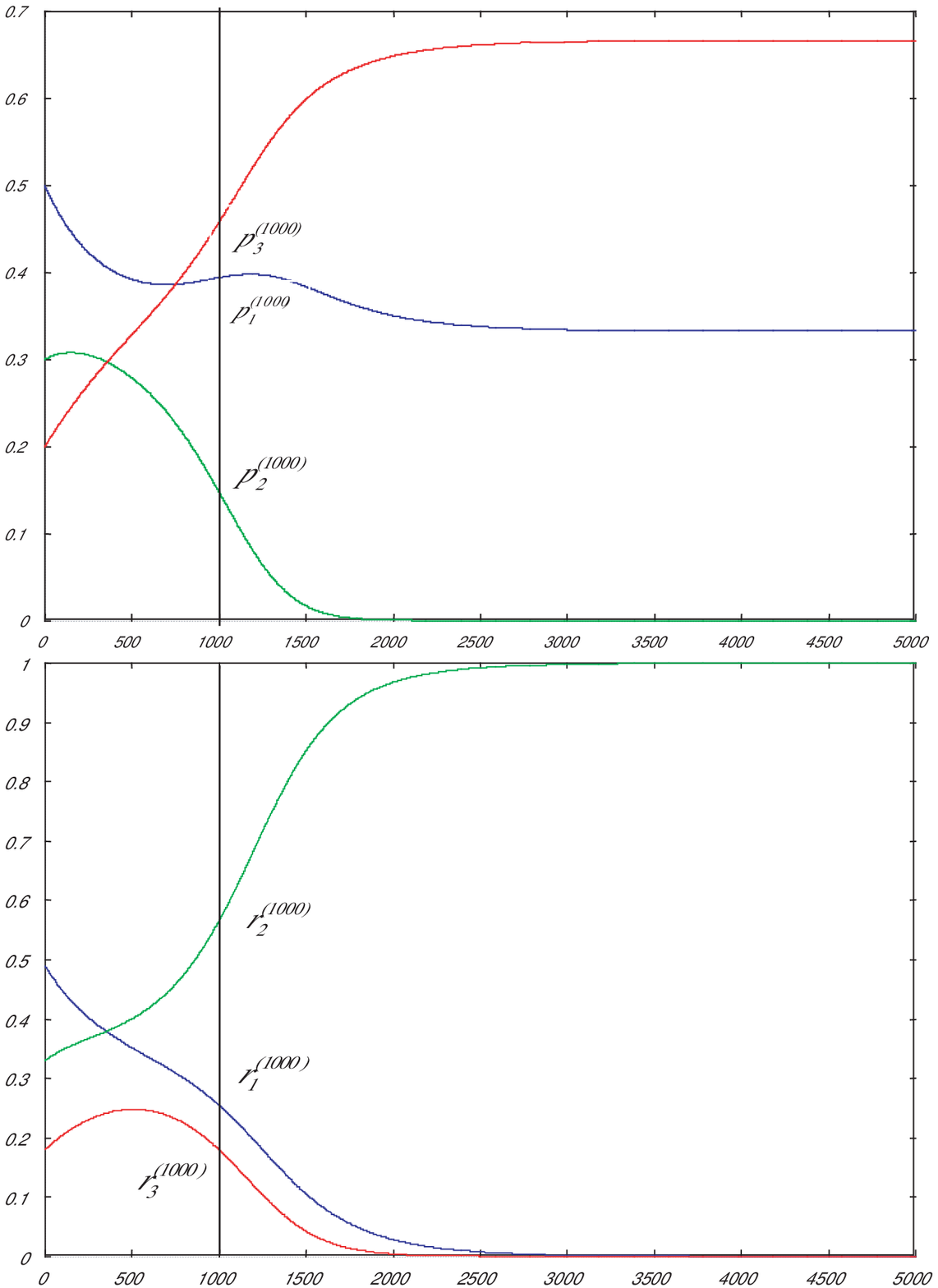}
\caption{$\alpha=1, \mathbf{p}^0=(0.5;0.3;0.2),
\mathbf{p}^0=(0.48;0.34;0.18)$} $\mathbf{p}^{\infty}=(0.33; 0;
0.67), \mathbf{r}^{\infty}=(0; 1; 0).$
\end{figure}

If we start with a pair of identical vectors,
$\mathbf{p}=\mathbf{r}$, then
$\mathbf{p}^{\infty}=\mathbf{r}^{\infty}$ too. That is, all
non-zero coordinates of the limiting vectors are equal.

In the general case, $\mathbf{p},\mathbf{r}\in \mathbb{R}^N_+,$
the coordinates $p_i^{(n)}, r_i^{(n)}$ have at most several
oscillations and then reach monotonically their positive or zero
limits. The limiting values $p_i^{\infty}, r_i^{\infty}$ may be
described in terms of starting states.

Given a couple of stochastic vectors
 ${\bf p}, {\bf r} \in {\mathbb R}_+^{n}, \ {\bf p}\neq {\bf r}, \ ({\bf p}, {\bf r})>0$,
 define
 $$D_+:=\sum_{i\in {\mathbb N}_+ }d_i, \
D_-:=\sum_{i\in {\mathbb N}_- }d_i, $$
 where
 $$ d_i=p_i-r_i, \ \ {\mathbb  N}_+:=\{i: d_i >0 \},\ {\mathbb  N}_-:=\{i: d_i <0 \}.$$

 Obviously $$0<D_+=-D_-<1,$$ since ${\bf p}\neq{\bf r}$, and
$\sum_ip_i-\sum_ir_i=0=D_++D_-$.

\ \ \ \ \ \ \ \ \ \ \

  {\bf Theorem 2.} {\it  Let ${\bf p}\neq {\bf r}, \ ({\bf p}, {\bf r})>0$.
   In the purely repulsive case,  $\alpha=1$,
  the
  coordinates of the limiting vectors
  ${\bf p}^{\infty},  {\bf
r}^{\infty}$ have the following explicit distributions:
\begin{equation}\label{result_p}
p_i^{\infty}=\left\{
\begin{array}{ll}
d_i/D, & \quad i\in {\mathbb N}_+ \\ 0, & \text{otherwise}
\end{array}\right., \ \ \
r_i^{\infty}=\left\{
\begin{array}{ll}
-d_i/D, & \quad i\in {\mathbb N}_- \\ 0, & \text{otherwise,}
\end{array}\right.
\end{equation}
where $D:=D_+=-D_-$.}

{\bf Remark.} From  (\ref{result_p}) it follows that any
transformation ${\bf p, r } \to {\bf p', r'}$, which does not
change the values $d_i$ and $D$, preserves the same limiting
distribution as for the  vectors ${\bf p}^{\infty}, {\bf r
}^{\infty}$. A class of such transformations may be presented by a
shift transformation of coordinates, $p_i \to p'_i=p_i+a_i, \ r_i
\to r'_i=r_i+a_i$ with  appropriated $a_i's$.

In the case $-1\leq\alpha<0$ of the pure attractive interaction we
have another limiting distribution.

Define the set $\mathbb{S}_0:=\{k|p_k^{\infty}=r_k^{\infty}=0\}$
and set $$\mathbb{S}^{\infty}:=\{1,\ldots,N\}\backslash
\mathbb{S}_0.$$

{\bf Theorem 3.} {\it In the purely attractive  case, $\alpha=-1$,
  the limiting vectors
  ${\bf p}^{\infty},  {\bf
r}^{\infty}$ are equal and their coordinates  have the following
distributions:
\begin{equation}\label{result_pa}
p_i^{\infty}=r_i^{\infty}=\left\{
\begin{array}{ll}
1/m, & \quad i\in {\mathbb S}^\infty \\ 0, & \text{otherwise,}
\end{array}\right. \ \ \
\end{equation}
where $m=|{\mathbb S}^\infty|$ denotes the cardinality of the set
${\mathbb S}^\infty$.}

In general, it is an open question to give a complete
characterization of $\mathbb{S}_0$.

Below we present several sufficient conditions for $k$ to belong
to the set $ {\mathbb S}_0$. Simultaneously these conditions give
some characterization  for the points to be in
 ${\mathbb S}^\infty$.

We will use the following notations:
\begin{equation}\label{result_par} \sigma_i:=p_i+r_i, \
\rho_i:=p_ir_i, \ \sigma_i^1:=p_i^{1}+r_i^{1} \
\rho_i^{1}:=p_i^{1}r_i^{1}.\end{equation}

\ \ \ \ \ \ \ \ \ \

 {\bf Proposition 1}. {\it If \begin{equation}\label{result_pan} \sigma_i\geq\sigma_k, \ \
\rho_i>\rho_k, \ \  {\rm or} \ \  \sigma_i>\sigma_k, \ \
\rho_i\geq\rho_k,\end{equation} then
$$p_k^{\infty}=r_k^{\infty}=0,$$ and therefore $k\in
\mathbb{S}_0.$}

\textit{Proof.} By (\ref{result_par}) we have
$$\sigma^1_k=p^1_k+r^1_k-1/z(p_k+r_k+2p_kr_k)=1/z(\sigma_k+2\rho_k)$$
where we recall that $z=1+(\mathbf{p},\mathbf{r}).$ Therefore each
of the conditions (\ref{result_pan}) implies that
$\sigma_i^1>\sigma^1_k$. Further, since
\begin{equation}\label{result_pal}\rho^1_k=1/z^2(\rho_k+(\rho_k)^2+\rho_k\sigma_k),\end{equation} again from
(\ref{result_pan}) it also follows that $\rho^1_i>\rho_k^1.$ Thus,
by induction, $\sigma_i^N>\sigma_k^N$ and $\rho_i^N>\rho_k^N$ for
all $N\geq 1.$

Or, in other words,
$$1<\frac{p_i}{p_k}<\frac{p_i^1}{p_k^1}<\ldots<\frac{p_i^N}{p_k^N}\ldots,$$
\begin{equation}\label{result_paq} 1<\frac{p_i}{p_k}<
\frac{r_i^1}{r_k^1}<\ldots<\frac{r_i^N}{r_k^N}\ldots, N=1,2,\ldots
\end{equation}

Thus, sequences of the ratios
$$\frac{p^N_i}{p^N_k}, \frac{r^N_i}{r^N_k}$$ are monotone
increasing as $N\to\infty.$ Assume for a moment that there exists
a finite limit,
$$1<\lim\limits_{N\to\infty}\frac{p^N_i}{p^N_k}=\frac{p^{\infty}_i}{p^{\infty}_k}\equiv
\frac{p^{\infty}_i}{p^{\infty}_k}\cdot\frac{1+r^{\infty}_i}{1+r^{\infty}_k}=M<\infty.$$
This is only possible if $r^{\infty}_i=r^{\infty}_k,$ which
contradicts (\ref{result_paq}). Thus, $M=\infty$ and therefore
$p_k^{\infty}=0$, as well as $r_k^{\infty}=0.$ $\Box$

Let us  consider now the critical situation, when for a fixed pair
of indices, say $i$ and $k$, the values $\sigma_k-\sigma_i$, \
$\rho_k-\rho_i$ have opposite signs, for example,
$\sigma_k-\sigma_i>0$, \ $\rho_k-\rho_i<0$. In such a case it is
not clear what behavior the coordinates $p_i^{N}, r_i^{N}$ and
 $p_k^{N}, r_k^{N}$ will have when
  $N\to \infty$. We will show that the
limits  depend on which of the two  values,
 $2\rho_i+\sigma_i$ or $2\rho_k+\sigma_k$, is larger.
Moreover we will show that even if  $p_k$ is the largest
coordinate, it may happen that $p_k^{\infty}=0$.
 Let for example,  $p_k=\max_j\{p_j,r_j\}$ and
 $\sigma_k=p_k+r_k>p_i+r_i=\sigma_i, $ however
 the value of $r_k$ is such  that $
\rho_k=p_kr_k < p_ir_i=\rho_i.$ \ Then  under some additional
condition it is possible to have $p_k^{\infty}=0$. In fact we
have:

\ \ \ \ \ \ \ \ \

 {\bf Proposition 2. } {\it Let for the  coordinates $p_i, r_i, p_k, r_k, \ i\neq
 k$,
  the following
  conditions be fulfilled:
\begin{equation}\label{c1}\sigma_k>\sigma_i\end{equation} but
\begin{equation}\label{c2}\rho_k < \rho_i.\end{equation} Assume
\begin{equation}\label{c3}
2\rho_k + \sigma_k \leq 2\rho_i+\sigma_i.
\end{equation}
Then
\begin{equation}\label{r3}
p_k^{\infty}=r_k^{\infty}=0,
\end{equation} i.e., $k\in \mathbb{S}_0$
}

\textit{Proof.} We will show that (\ref{c1}), (\ref{c2}), and
(\ref{c3}) imply \begin{equation}\label{r4}
p_k^1+r_k^1=\sigma_k^1\leq\sigma_i^1=p_i^1+r_i^1
\end{equation}
and \begin{equation}\label{r5}
p_k^1r_k^1=\rho_k^1<\rho_i^1=p_i^1r_i^1.
\end{equation}
Then (\ref{r3}) follows from Proposition 1. In reality (\ref{r4})
follows from (\ref{c3}) directly, without condition (\ref{c2}).
So, we have only to prove (\ref{r5}).

With this aim we find the representation for $\rho_i^1$ in terms
$\sigma_i$ and $\sigma_i^1$. Since
$\sigma_i^1=1/z(\sigma_i+2\rho_i)$ we have
\begin{equation}\label{fff}
\rho_i=1/2(z\sigma_i^1-\sigma_i).
\end{equation} By (\ref{result_pal}) and (\ref{fff}) we get $$\rho_i^1=1/z^2(\rho_i+\rho_i^2+
\rho_i\sigma_i)=\frac{1}{2z^2}(z\sigma_i^1-\sigma_i)[1+1/2(z\sigma_i^1-\sigma_i)+\sigma_i]$$
$$=\frac{1}{4z^2}(z\sigma_i^1-\sigma_i)(2+z\sigma_i^1+\sigma_i)=\frac{1}{4z^2}[2z\sigma_i^1+z^2(\sigma_i^1)^2+z\sigma_i^1\sigma_i-2\sigma_i-z\sigma_i^1\sigma_i-\sigma_i^1]$$
$$=\frac{1}{4z^2}[2z\sigma_i^1+z^2(\sigma_i^1)^2-\sigma_i^2-2\sigma_i].$$

Therefore
$$\rho_k^1-\rho_i^1=1/z^2[\rho_k(1+\rho_k+\sigma_k)-\rho_i(1+\rho_i+\sigma_i)].$$

Thus, we have
$$\rho_k^1-\rho_i^1=1/4z^2[2z(\sigma_k^1-\sigma_i^1)+z^2((\sigma_k^1)^2-(\sigma_i^1)^2)+((\sigma_i)^2-(\sigma_k)^2)+2(\sigma_i-\sigma_k)]<0$$
due to starting condition (\ref{c2}), and (\ref{r4}). Thus
$\rho_k^1<\rho_i^1,$ i.e., (\ref{r5}) is true. $\Box$

 We
stress that  (\ref{r3}) is true in spite of $\sigma_k>\sigma_i$.
Of course, if $\sigma_k<\sigma_i$ and $\rho_k<\rho_i$, then
(\ref{r3})  holds without any additional condition of the form
(\ref{c3}).

\section{Model of conflict interaction between complex systems}
In this section we construct a dynamical model of conflict
interaction between a pair of complex systems. Each of the systems
is subjected to the inner conflict between their elements. For
simplicity, we assume both systems to be similar and described by
discrete prey-predator models of type (2.5). We introduce the
conflict interaction between these systems using an approach
developed in \cite{ABodK, AKPT2, BKoDo, KoTC, KoTC1, KoDo, KoKh}.
With such a rather complex situation we may obtain a wide spectrum
of evolutions. In this work we study qualitative characteristics
of the behavior of corresponding dynamical systems for some choice
of parameters $a,b,c,d,e,f,\alpha$ (see (2.5), (3.1)) and values
of initial populations of species $P_i, R_i$.

The coefficient $\alpha$, that shows intensity of the interaction
between systems, has an important effect. The increasing $\alpha$
from zero to unit causes the appearance of a series of bifurcations.
For $\alpha=0$ we have two copies of independent Lotka-Volterra
models. For small values of $\alpha$ both systems behave like pure
Lotka-Volterra systems, coming them to a stable state.


Under fixed parameters and the starting coordinates $a=0.2, b=0.006,
c=0.002, d=0.008, e=0.002, f=0, P_1^{(0)}=3, P_2^{(0)}=10,
R_1^{(0)}=5, R_2^{(0)}=20$
we have first bifurcation
point at $\alpha\approx 0.0056781739$. The coordinates
$P_i^{(n)}(R_i^{(n)})$ oscillate and a cycle of a small period
appears.


The following increase of $\alpha$ shows the appearance of new
bifurcation points that are characterized by an increasing value
of the cyclic period. For the value $\alpha=0.4815545975$ a cycle
of infinite period appears. This means that all coordinates
rapidly reach the stable state. In this case some species may
disappear, even if they had some stable positive values in a pure
($\alpha=0$) Lotka-Volterra model.

The role of the coefficients $a,b,c,d,e,f$ and initial quantity of
the species $P_i, R_i$ in a pure Lotka-Volterra model is
well-known and described (see, e.g., \cite{Mur, KuBe}). Partially,
coefficients $a,d$ govern the increase of the pray population when
predators are absent and the predator population decreasing when
prays are absent. In turn, the coefficients $b,e$ are responsible
respectively for the pray quantity decreasing with an increasing
number of predators, and increase of the predator population with
an increase of the number of prays. The last coefficients in each
of the equations give the limitation of increasing of both
populations. In other words, each population "makes pressure" on
itself, it does not permit an infinite reproducibility.

Questions about stable points, orbits, asymptotic behavior of
orbits are well described for the classical Lotka-Volterra model.
We shall recall that usually there are at least three equilibrium
points. They are refereed in literature as follows (see, e.g.
\cite{Mur}):

(1) trivial (0,0);

(2) axial ($a/b$,0);

(3) inner positive
$$\left(\frac{a}{b}-\frac{b}{c}\frac{ae-cd}{be+cf},\frac{ae-cd}{be+cf}\right).
\eqno(4.1)$$

An equilibrium point is called stable point if after a sudden
change of population it comes back to an equilibrium point some
time later. This may happen monotonically, or with some
oscillations.

We should note that under the existence of stable points the
behavior of the system is well defined by coefficients
$a,b,c,d,e,f$. But under the absence of stable points, the
behavior of the system is defined by the initial data $P_i, R_i$.
Depending on how close the initial data are situated with respect
to the equilibrium point, the system may evaluate in a different
way.

The role of all these coefficients is preserved in the case of our
model. But now their influence is much more complex. We present
here only first steps in this direction. We shall discuss not only
stability zones, as it was pointed above, but also the values of
the coefficients for which the system oscillates along some closed
cycles.

The state of our dynamical system is fixed by a pair of vectors
$\mathbf{P}^n=(P_1^{(n)}, \ldots, P_N^{(n)}),$
$\mathbf{R}^n=(R_1^{(n)}, \ldots, R_N^{(n)})$ with non-negative
coefficients, where $n=0,1,\ldots$ denotes the discrete time,
$N\geq 2$ stands for the number of conflict positions. Here we
study the most simple situation, when every system consists of
only two agents: pray and predator, i.e. $N=2$. The complex
conflict transformation is denoted by the mapping
$$\begin{pmatrix}
  \mathbf{P}^{n} \\
  \mathbf{R}^{n}
\end{pmatrix}\begin{array}{c}
  \stackrel{F}{\longrightarrow}
\end{array}\begin{pmatrix}
  \mathbf{P}^{n+1} \\
  \mathbf{R}^{n+1}
\end{pmatrix},$$ where $F$ is the composition of four operations, the specific mathematical
transformations: $F=[\mathcal{N}^{-1}*\mathcal{N}]U.$

Let us describe them in an explicit form for the first step.

The first operation $U$ describes the interaction between elements
inside every system separately according to the pray-predator
model. Corresponding mathematical transformation of vectors (the
interaction composition)
$\{\mathbf{P}^0,\mathbf{R}^0\}\stackrel{U}{\longrightarrow}\{\mathbf{\widetilde{P}}^0,\mathbf{\widetilde{R}}^0\}$
is described by the system of equations of the form (2.5):
 $$\begin{array}{c}
  \widetilde{P}_1^{(0)}=P_1^{(0)}+P_1^{(0)}(a-bP_2^{(0)}-cP_1^{(0)}), \\
  \widetilde{P}_2^{(0)}=P_2^{(0)}+P_2^{(0)}(-d+eP_1^{(0)}-fP_2^{(0)}),
\end{array}$$ and $$\begin{array}{c}
  \widetilde{R}_1^{(0)}=R_1^{(0)}+R_1^{(0)}(a-bR_2^{(0)}-cR_1^{(0)}), \\
  \widetilde{R}_2^{(0)}=R_2^{(0)}+R_2^{(0)}(-d+eR_1^{(0)}-fR_2^{(0)}),
\end{array}$$ where the passage to new values of coordinates is pointed by tilde,
but not by changing of upper index, likely to (2.5).

The following operation involves the interaction $*$ (see (3.1))
between previous systems according to the theory of the
alternative conflict for non-annihilating opponents (see, e.g.
\cite{ABodK, AKPT2, BKoDo, KoTC, KoTC1, KoDo, KoKh}). To describe
this operation we at first have to normalize the vectors
$\mathbf{\widetilde{P}}^0=(\widetilde{P}_1^{(0)},
\widetilde{P}_2^{(0)}),
\mathbf{\widetilde{R}^0}=(\widetilde{R}_1^{(0)},
\widetilde{R}_2^{(0)})$, i.e., to work with stochastic vectors.

We use the following notation for normalization:
$\mathcal{N}\{\mathbf{\widetilde{P}}^0,\mathbf{\widetilde{R}}^0\}=
\{\mathbf{p}^0,\mathbf{r}^0\},$ where the coordinates of the
stochastic vectors $\mathbf{p}^0,\mathbf{r}^0$ are determined by
formulae $${p}_1^{(0)}=
\frac{\widetilde{P}_1^{(0)}}{\widetilde{z}^{(0)}_P},
{p}_2^{(0)}=\frac{\widetilde{P}_2^{(0)}}{\widetilde{z}^{(0)}_P},
{r}_1^{(0)}=\frac{\widetilde{R}_1^{(0)}}{\widetilde{z}^{(0)}_R},
{r}_2^{(0)}=\frac{\widetilde{R}_2^{(0)}}{\widetilde{z}^{(0)}_R},$$
where $\widetilde{z}^{(0)}_P=\widetilde{P}_1^{(0)}+
\widetilde{P}_2^{(0)},
\widetilde{z}^{(0)}_R=\widetilde{R}_1^{(0)}+
\widetilde{R}_2^{(0)}.$

The next step exactly corresponds to the conflict interaction
between systems. We introduce new stochastic vectors
$\{\mathbf{p}^1, \mathbf{r}^1\}$ with coordinates:
$${p}_j^{(1)}=\frac{{{p}}_j^{(0)}(1-\alpha
{{r}}_{j}^{(0)})}{1-\alpha\sum_{i=1}^2{{p}}_i^{(0)}{{r}}_i^{(0)}},
{r}_j^{(1)}=\frac{{{r}}_j^{(0)}(1-\alpha
{{p}}_{j}^{(0)})}{1-\alpha\sum_{i=1}^2{{p}}_i^{(0)}{{r}}_i^{(0)}},
j=1,2. $$

Finally, we have to come back to the non-normalized vectors, which
characterize quantitatively populations in both regions after
inner and outer conflicts operations. So, at time $n=1$ we have
the following vectors
$\mathcal{N}^{-1}\{\mathbf{p}^1,\mathbf{r}^1\}=\{\mathbf{{P}}^1,\mathbf{R}^1\}$,
where $$\mathbf{P}^1=(P_1^{(1)},
P_2^{(1)}),\mathbf{R}^1=(R_1^{(1)}, R_2^{(1)}),$$ and where
$${P}_j^{(1)}={p}_j^{(1)}\widetilde{z}^{(0)}_P,
{R}_j^{(1)}={r}_j^{(1)}\widetilde{z}^{(0)}_R, j=1,2.$$

We can repeat this procedure starting from
$\{\mathbf{{P}}^1,\mathbf{R}^1\}$. So we get
$\{\mathbf{{P}}^2,\mathbf{R}^2\}$. And so on for any $n$th step.

To find the equilibrium points in the case of the complex conflict
interaction described above, we have to solve the following system
of equations for $P_1, P_2, R_1, R_2$: $$\left\{\begin{array}{c}
 (a+1-bR_1-cP_1)(Z_2-\alpha R_2(-d+1+eP_2-fR_2))Z_1=Z, \\
 (-d+1+eP_1-fR_1)(Z_2-\alpha P_2(a+1-bR_2-cP_2))Z_1=Z, \\
 (a+1-bR_2-cP_2)(Z_1-\alpha R_1(-d+1+eP_1-fR_1))Z_2=Z, \\
 (-d+1+eP_2-fR_2)(Z_1-\alpha P_1(a+1-bR_1-cP_1))Z_2=Z, \\
  \end{array}\right.$$ where
$$Z_1=P_1(a+1-bR_1-cP_1)+R_1(-d+1+eP_1-fR_1),$$
$$Z_2=P_2(a+1-bR_2-cP_2)+R_2(-d+1+eP_2-fR_2),$$
$$Z=Z_1Z_2-\alpha[P_1P_2(a+1-bR_2-cP_2)(a+1-bR_1-cP_1)+$$ $$+
R_1R_2(-d+1+eP_1-fR_1)(-d+1+eP_2-fR_2)].$$

We note that in the case $\alpha=0$ we have two copies of pure
Lotka-Volterra models and the corresponding system of equations
has at least three equilibrium points (trivial, axial, inner
positive).

 For the case $\alpha\neq0$ it is difficult
 to obtain exact solutions. Let us obtain some insights by numerical approximation.

 Partially, we found that there exist equilibrium points
 and the limit cycles for a wide set of parameter values and initial data (see Figure  4-7).

\begin{figure}
\centering
\includegraphics{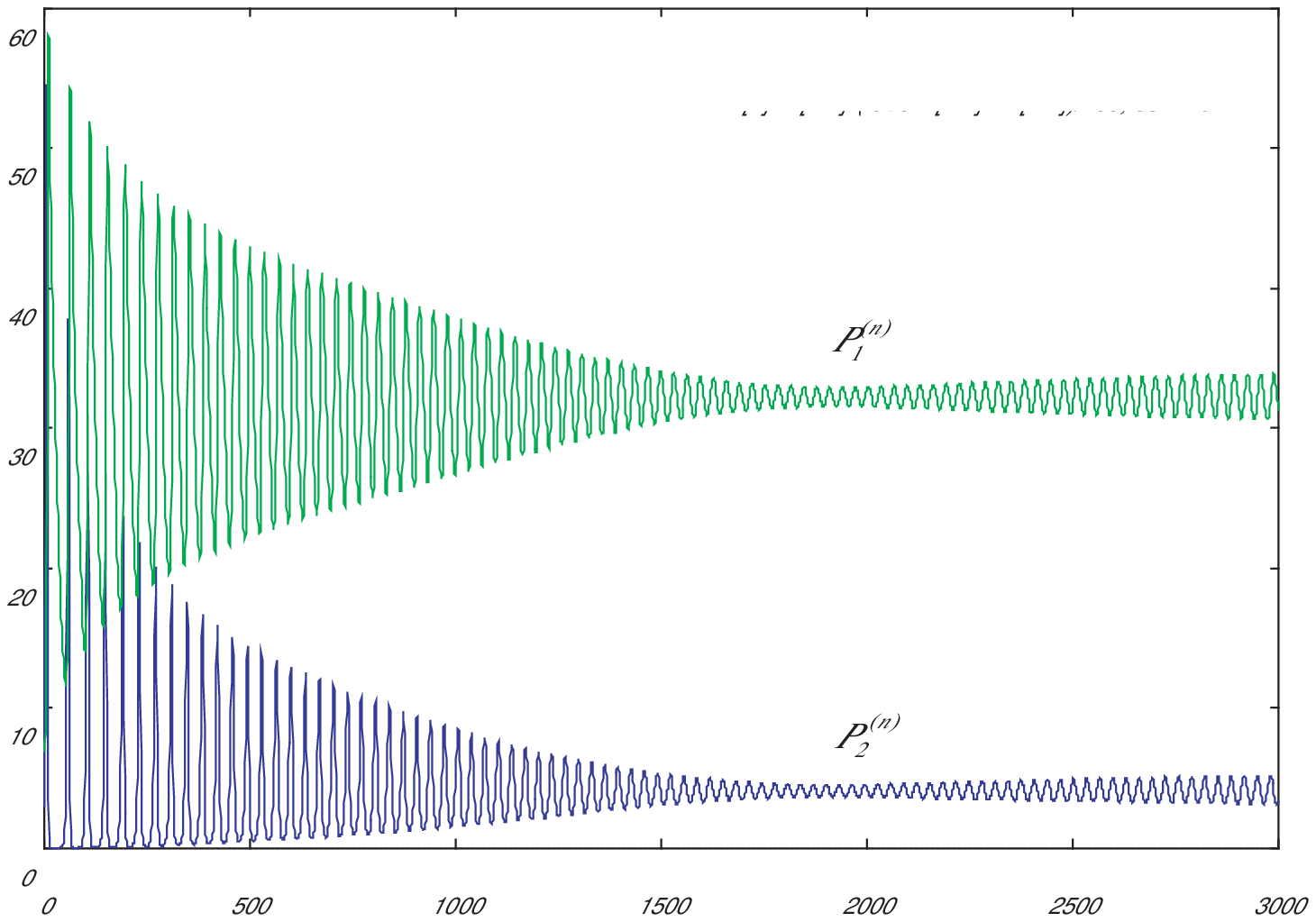}
\caption{{\it The existence of the strong bifurcation produces
oscillations of the large}} {\it amplitude. $a=0.2, b=0.006,
c=0.002, d=0.008, e=0.002, f=0,$ the conflict interaction
coefficient} $\alpha=0.01,$ $$P_1^{(0)}=3, P_2^{(0)}=5,$$
$$R_1^{(0)}=7, R_2^{(0)}=10.$$
\end{figure}

\begin{figure}
\centering
\includegraphics{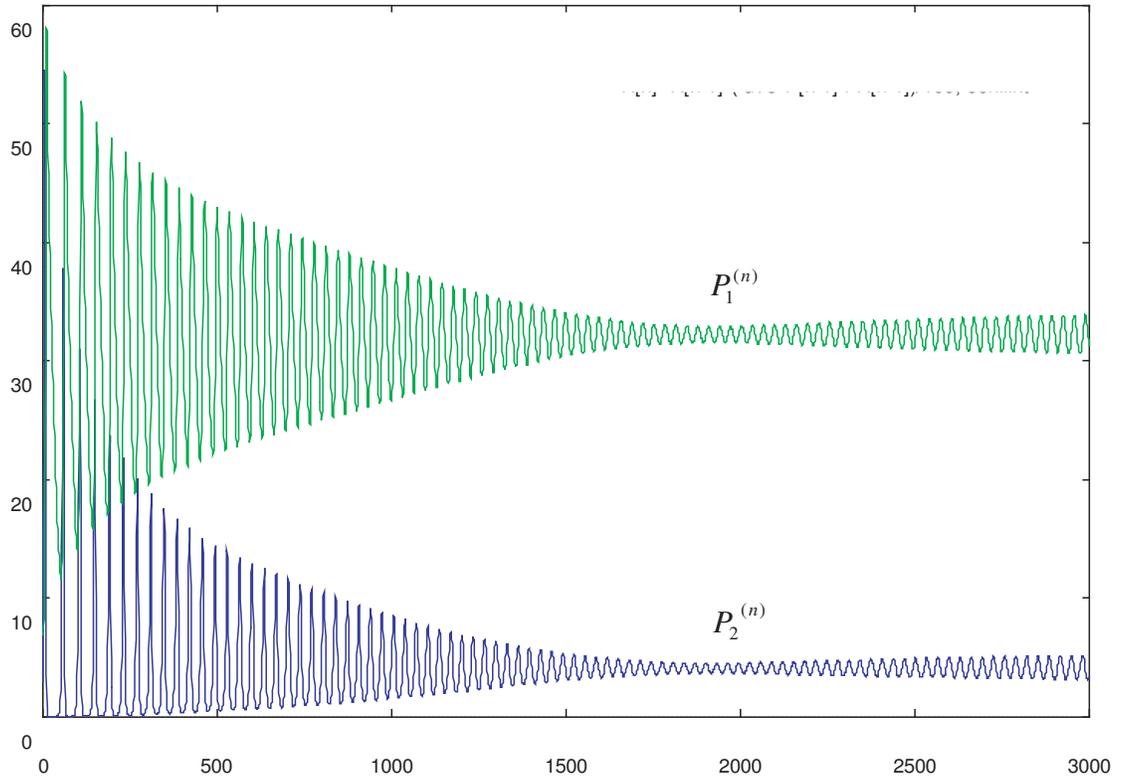}
\caption{{\it The corresponding phase-space $(P_1^{(0)},
P_2^{(0)})$.}}
\end{figure}

\begin{figure}
\centering
\includegraphics{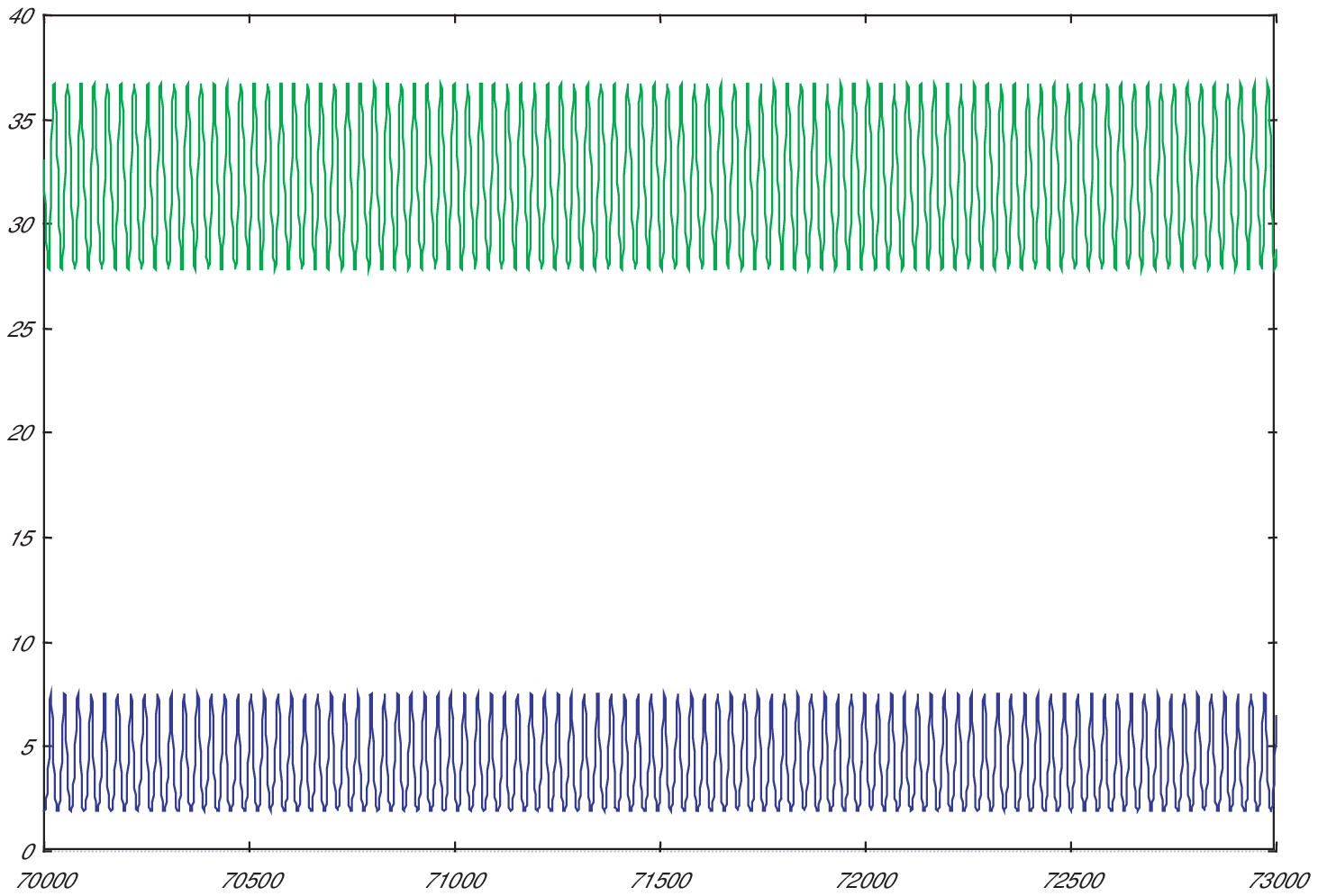}
\caption{{\it The existence of the stable oscillations of the
conflict interaction}} {\it between Lotka-Volterra systems after
70000 steps of iteration.}
\end{figure}

\begin{figure}
\centering
\includegraphics{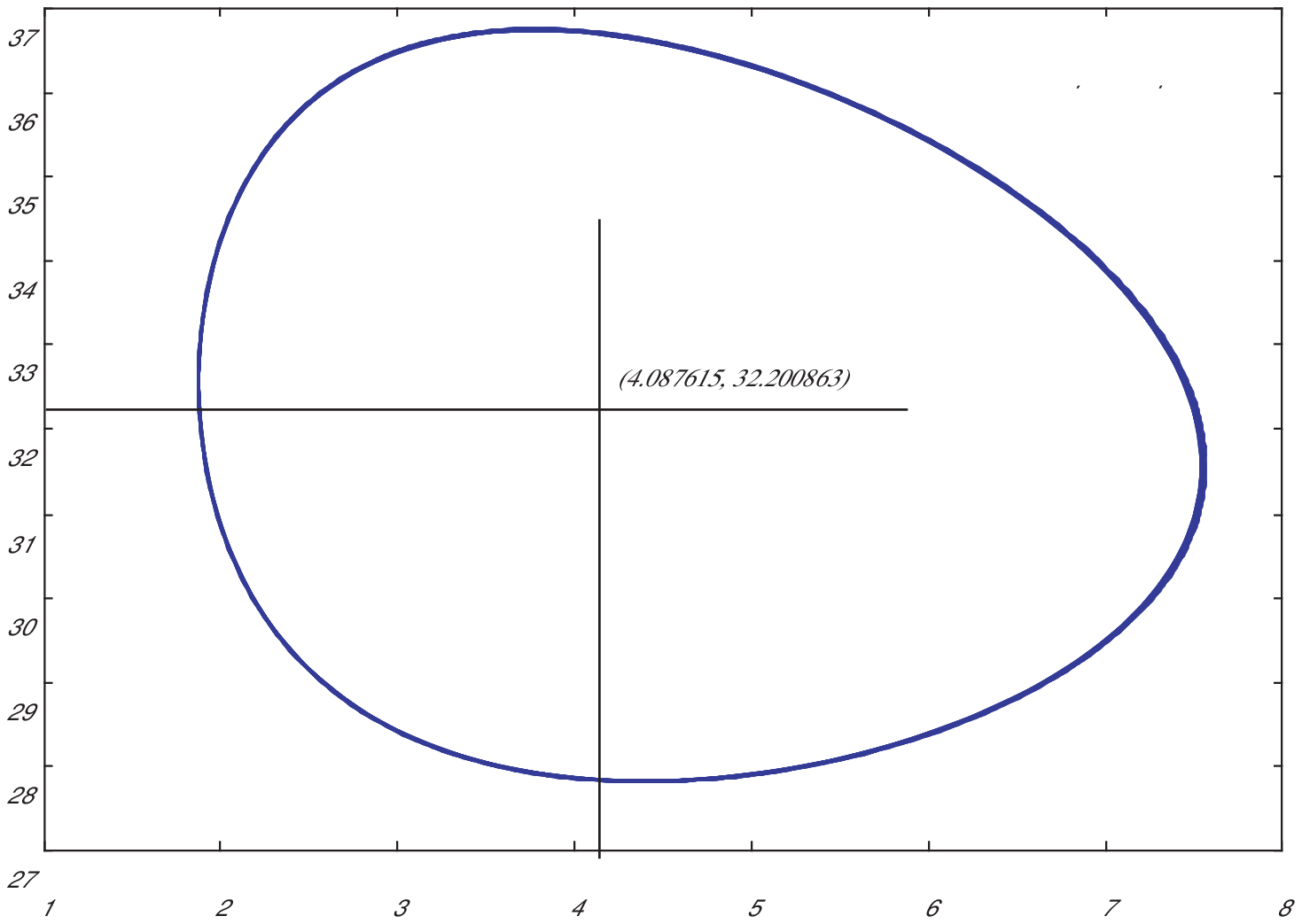}
\caption{{\it The limiting cycle in the corresponding phase-space
$(P_1^{(0)}, P_2^{(0)})$ after}} {\it 70000 steps of iteration.
The unstable equilibrium point is shown.}
\end{figure}

Moreover, we established the shift effect for the equilibrium point.
Namely, we observe that the inner positive equilibrium point (it
exist in any system and may be found by formula (4.1)) is shifted
after the application of the conflict interaction between systems.
We see by (4.1) that stabilization of discrete Lotka-Volterra model
with parameters $a=0.2, b=0.006, c=0.002, d=0.008, e=0.002, f=0,
\alpha=0.007, P_1^{(0)}=3, P_2^{(0)}=10, R_1^{(0)}=5, R_2^{(0)}=20$
occurs when $P_1=4, P_2=32.$ This may be easily verified by putting
these initial data into corresponding equations. In this case we
have trivial dynamics.

Let us consider the case of discrete Lotka-Volterra model with the
conflict interaction between systems. We take the values of the
coefficients $a=0.2, b=0.006, c=0.002, d=0.008, e=0.002, f=0,
\alpha=0.005$. Now the equilibrium point has the coordinates
$P_1=4.043507, P_2=32.100629.$
 The dynamics is constant with these initial data.

In case of larger $\alpha$, when oscillations appear, the
equilibrium point may also be easily found if we put the initial
data in both
 systems to be equal. In this case the behavior is
 like in the case of a pure Lotka-Volterra model, and stabilization occurs. However, the
 stable point is shifted, for example, when $\alpha=0.01$ (see Figure 4-7) the equilibrium point is
 $P_1=R_1=4.087615, P_2=R_2=32.200863.$

 Thus, if we have some pray-predator system and want to change the population inside this system,
  we may create an analogous "artificial" system, introduce the conflict interaction and obtain
  the desired
  shift of the equilibrium point. Apparently a stronger shift
  of the stable equilibrium point in the case of an "ensemble" of larger
  amount of Lotka-Volterra systems. So, we observed the interesting phenomenon:
  the equilibrium point of an isolated system is shifted if we
  come to the case when
  identical systems are united as an "ensemble".

 However, this equilibrium point is unstable, any perturbation of initial data causes the receding of the system from
 the equilibrium point.

 One of more interesting observations concerns the limit cycles.
 It is known that no such kind of orbits in discrete Lotka-Volterra model is possible.
  But under the effect of the outer conflict, as we see at the pictures, the dynamical system reaches
the limit cycle starting both from an inside or outside point with
respect to the orbit.
    Partially, in Figure 10,11 we present the model, that starts at $P^{(0)}_1=4, P^{(0)}_2=32$.
    As it was pointed above, in case of a pure Lotka-Volterra model, with these initial data there is
    no dynamics.
However, in the case of the model with the outer conflict the
process tends to a limit cycle.

\begin{figure}
\centering
\includegraphics{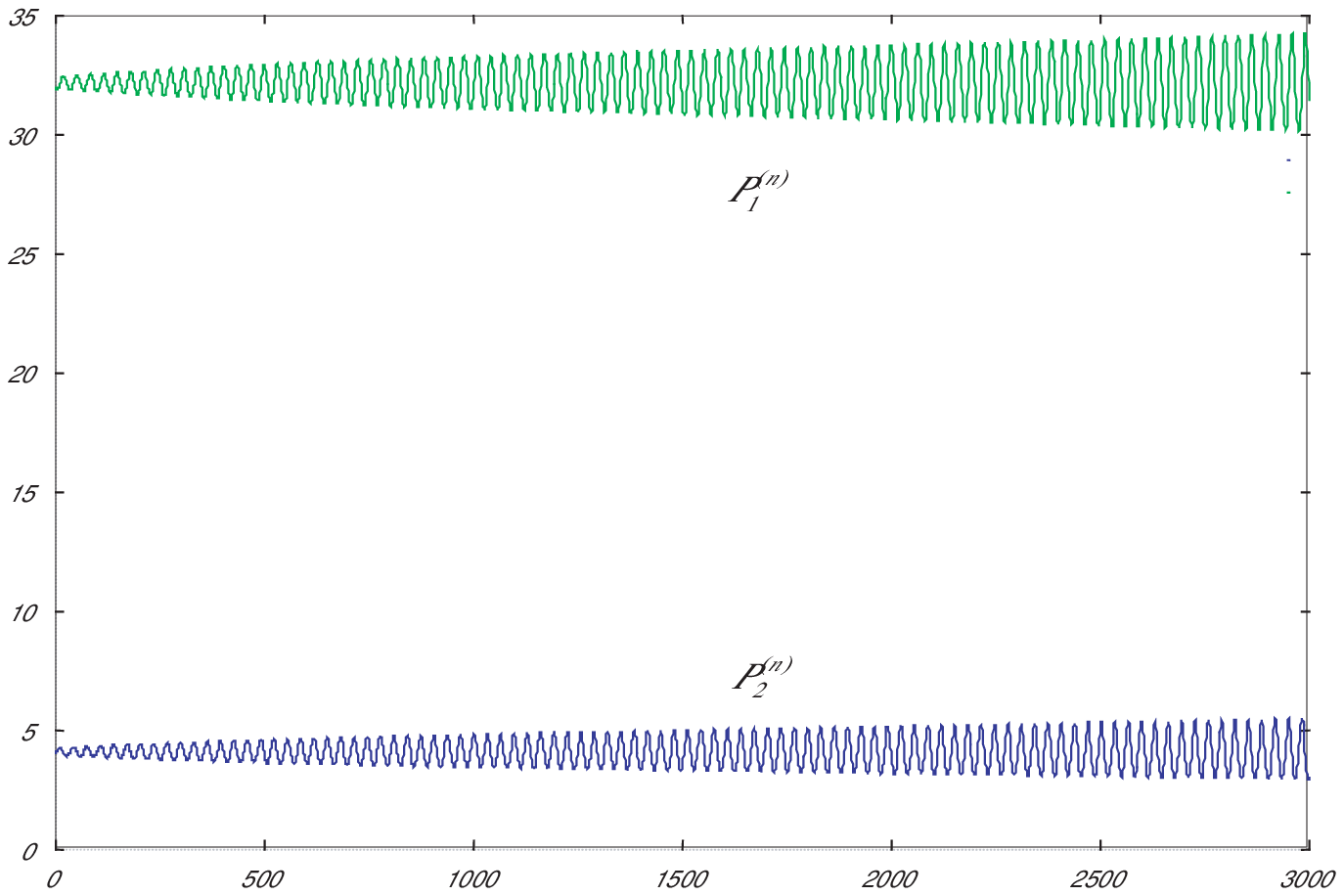}
\caption{{\it Conflict interaction between Lotka-Volterra systems.
The starting }} {\it parameters are the same as in Figure 4.}

{\it Initial data are inside the limit cycle (Figure 7).}
\end{figure}

\begin{figure}
\centering
\includegraphics{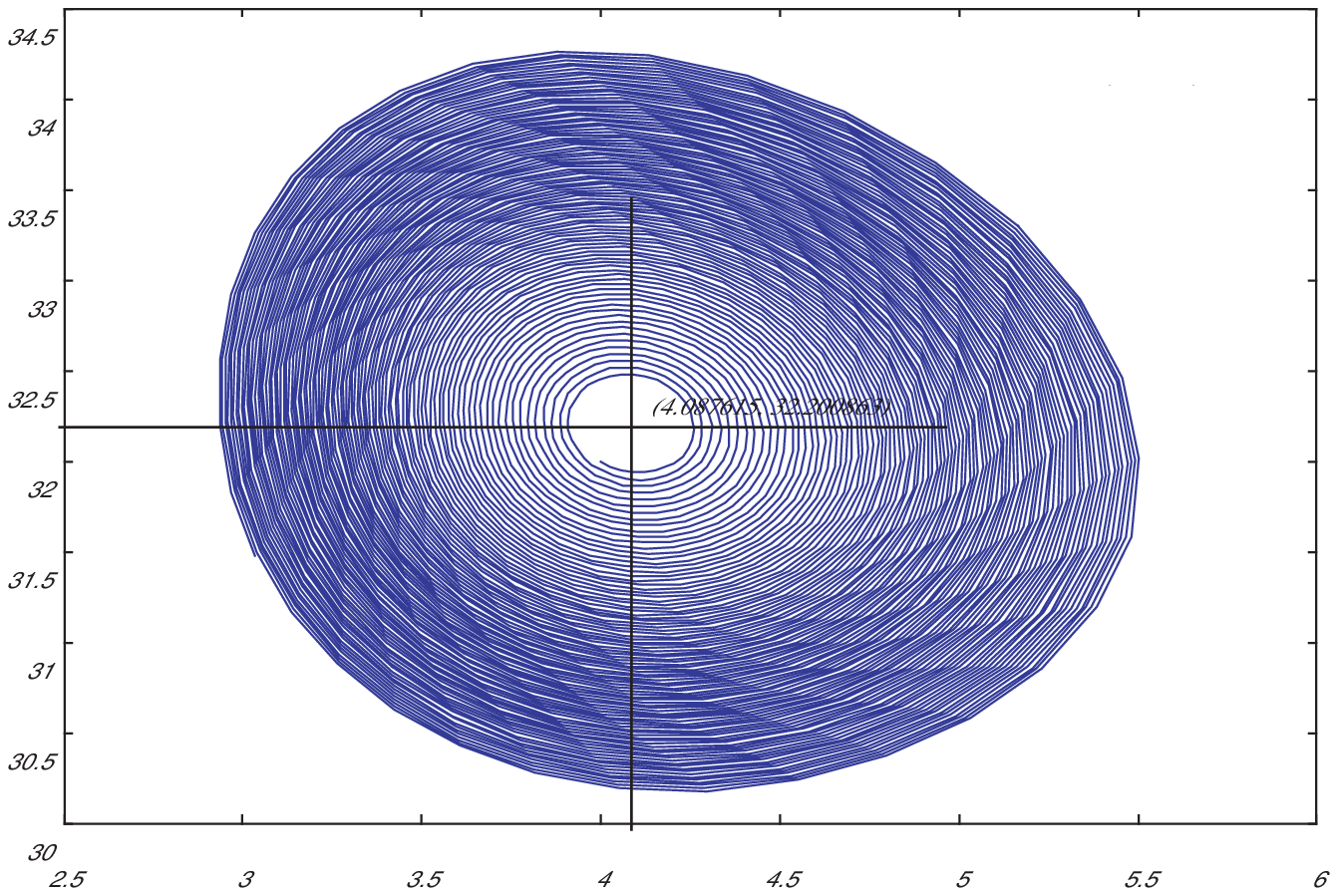}
\caption{{\it The corresponding phase-space $(P_1^{(0)},
P_2^{(0)})$. Trajectory tends to the }} {\it limit cycle, which is
an attractor. }

{\it Initial data are inside the limit cycle (Figure 7).}
\end{figure}

\section{Interpretation}
In many works on mathematical biology and economics \cite{BC, CM,
CKG, Ku, KuBe, LSGM, Mur, StOl, Tu} the modelling of population
dynamics or economical processes is based on Lotka-Volterra
equations. As a rule, continuous, not discrete, models are
studied. In some works the migration process is considered. It
takes place between different regions, inside which an interaction
of the Lotka-Volterra type is present. For example, in \cite{CM}
the migration rate between regions has some fixed probability.

We study discrete Lotka-Volterra models with an additional
interaction between them. That may be interpreted as a some kind
of correlation between the habitants of different regions. We
suppose that discrete models are more natural, partially it is
clear that birth and death of individuals happen at some fixed
moments of time.

It is well known that in the classical discrete pray-predator
model a stable point exist. The amount of prays and predators
tends to this point in the phase-space. In this case we observe
the following dynamics, after several period of oscillations the
populations stabilize (see Figure 2). Thus, we have an attracting
point in phase-space. Such a dynamics exists inside every region
when "migration" is absent.

When we introduce an additional interaction between the habitants
of different regions a redistribution process appears which we
interpret as a migration. In some of our complex models there is
no stable point, the amount of prays and predators in both regions
oscillates along fixed orbits.Appearently these orbits in a
phase-space are attractors.

We note that explicit formulas of conflict interaction between
non-annihilating opponents which describe the redistribution of
populations are given by (3.1). The individuals of a certain kind
migrate to the region, where their amount more numerous.

Is the "migration strategy" which is described in our model a
natural one? We suppose that in many cases individuals may be
right behaving in such a way. If we consider a pray-predator
model, it is clear that every separated individual is unable to
estimate all factors that have an influence on the population
dynamics like vital resources inside region, real amount of own
and alternative population, current population dynamics. In other
words, the individual "does not know" the parameters of the
Lotka-Volterra equations and their current influence on the
population dynamics.

However the individual has the group reflex and will migrate to
the region, where, as he supposes, the vital conditions are best
(his population should be concentrated there). He suggests, right
there are the resources, possibilities for reproduction, better
conditions to organize large groups. Formula (3.1) just describes
this tendency.

Similar motivations may be proposed in case of the work migration.
Here the unemployed may be regarded as playing the role of
"prays", employees as playing the role of "predators". People, who
seek for work and migrate to another country, do not know, as a
rule, the real situation in the opposite region. They prefer to
migrate to the country where the majority of their friends
migrated (group reflex).

Similar, but opposite picture happens with employees who inverts
their capital to the region with a higher profit.

So, at the cost of migration accelerates the increasing of one of
species' population in one of the regions. But at the same time
there is an effect of the inner pray-predator "fight" inside every
system. Partially, the population influences itself at the cost of
the last term in Lotka-Volterra equation.As a result, some time
later the backward migration starts.

In the Figure 6 we may see the effect of delay, when the amount of
prays inside the region decreases, but the predators continue
migration to this region, until their amount starts decreasing by
following the Lotka-Volterra model.

We emphasize, that in our model, in comparison with discrete
Lotka-Volterra model, a cyclic oscillations of populations are
observed. Moreover, a cyclic attractor exists in the phase-space,
and the pray-predator trajectory tends to this orbit both from
inside or outside point with respect to this cycle (Figure 8,9).

We remark that in our model the normalization was fulfilled by the
amount of habitants of the region, so the component of the
corresponding vector may be large both at the cost of large
population of fixed individuals and at the cost of small whole
population of the region. So, a migration to the region with a lot
of "free space" is also possible.

We also studied model with the attracting interaction
($\alpha<0$). In this case we obtained formally a similar
dynamics, but now individuals migrate to the region where they are
less numerous. Such a migration strategy might be also natural for
some species, e.g. for individuals who hunt separately, control
large territory and have confrontation with relatives.


\begin{thebibliography}{99}

\bibitem{ABodK} S. Albeverio, M.  Bodnarchyk, V. Koshmanenko,
 Dynamics of  Discrete Conflict Interactions
 Between Non-annihilating Opponents, {\it MFAT}, {\bf 11}, no.~4, (2005).

\bibitem{AKPT2}  S. Albeverio, V. Koshmanenko, M. Pratsiovytyi, G. Torbin, Spectral
properties of image measures under infinite conflict interaction,
 {\it MFAT}, N 2, (2006).

\bibitem{BC}  M.Bandyopadhyay, J.Chattopadhayay, Ratio-dependent predator-prey model:
 effects of environmental fluctuation and stability,
{\it Nonlinearity}, No. 18, 913--936, (2005).

\bibitem{BKoDo} M.V.  Bodnarchyk, V.D. Koshmanenko, N.V. Kharchenko,
Properties of limit states of dynamical conflict system, {\it
Nonlinear oscillations}, {\bf 7}, N 4, (2004) 446-461. (in
Ukrainian)

\bibitem{CM}  A.Colato, S.S.Mizrahi, Effects of random migration in population dynamics,
 {\it Physical Review E}, {\bf 64}, 1--14, (2001).

\bibitem{CKG}  R.Cressman, V.Krivan, J.Garay, Ideal free distributions, evolutionary
 games, and population dynamics in multiple-species environment,
 {\it The American Naturalist}, {\bf \164}, No. 4, 473--489, (2004).

\bibitem{KoTC}  V. Koshmanenko, On the Conflict Theorem for a Pair of
Stochastic Vectors, {\it Ukrainian Math. J.}, {\bf \ 55}, No. 4,
(2003).

\bibitem{KoTC1}  V. Koshmanenko, The Theorem of Conflict for Probability
Measures,  {\it Math. Methods of Operations Research} \ \ {\bf
59}, No.2 (2004) 303--313.

\bibitem{KoDo} V.D. Koshmanenko, N.V. Kharchenko,
Invariant points of dynamical conflict system in the space of
piecewise uniformly distributed measures,  {\it Ukrainian Math.
J.}, {\bf 56}, N 7,
 (2004) 927--938.

\bibitem{KoKh} V. Koshmanenko, N.Kharchenko,
Spectral properties of image measures after conflict interactions,
 {\it Theory of Stochastic Processes},   {\bf 10}(26), N 3-4, 73-81, (2004).

 \bibitem{Ku}  Y.Kuang, Basic properties of mathematical population models,
{\it J. Biomath}, No. 17, 129--142, (2002).

\bibitem{KuBe}  Y.Kuang, E.Beretta, Global qualitative analysis of a ratio-dependent
  predator-prey system, {\it J. Math. Biol.}, , No. 36, 389--406, (1998).

  \bibitem{LSGM}  Y.Lonzonn, S.Solomon, J.Goldenberg, D.Mazarsky, World-size
 global markets lead to economic instability, {\it Acrificial life}, 357--370, (2003).

\bibitem{Lot}  A.J. Lotka, Relation between birth rates and death
rates, {\it Science}, 26, 21-22, (1907).

 \bibitem{Mal}  T.R. Malthus, An essay on the principle of population,
 Reprinted by Macmillan, (1894).

 \bibitem{Mur}  J.D. Murray, Mathematical biology I:An Introduction, Springer, (2002).

 \bibitem{SaTa1}  K.MD.M.Salam, K.I.Takahashi, Mathematical model of conflict and cooperation with non-annihilating
multi-opponent, {\it J. Interdisciplinary Math.}, \textbf{in press}.

 \bibitem{StOl}  L.Stone, R.Olinky, Phenomena in ecological systems, Experimental Chaos:
6-th Experimental Chaos Conference, 476--487, (2003).

\bibitem{SaTa}  K.I.Takahashi, K.MD.M.Salam, Mathematical model of conflict with non-annihilating
multi-opponent, {\it J. Interdisciplinary Math.}, {\bf 9}, No. 3,
459-473, (2006).

\bibitem{Tu}  J.Tufto, Effects of releasing maladapted individuals: a
demographic-evolutionary model, {\it The American Naturalist},
{\bf \ 158}, No. 4, 331--340, (2001).

\bibitem{Ver} P.F. Verhulst, Notice sur la loi que la population suit dans son
accroissement, In Correspondence math\'{e}matique et physique
publi\'{e}e par A.Qu\^{e}telet, {\bf 10}, 113-121, (1838).

\bibitem{Vol}  V. Volterra, Sui tentativi di applicazione della matematiche alle scienze
biologiche e sociali, {\it Giornale degli Economisti}, {\bf 23},
436-458 (1901).

\end{thebibliography}
\end{document}